# The zeros of the Riemann-zeta function and the transition from pseudo-random to harmonic behavior

R. V. Ramos

rubens.viana@pq.cnpq.br

*Lab. of Quantum Information Technology, Department of Teleinformatic Engineering – Federal University of Ceara - DETI/UFC, C.P. 6007 – Campus do Pici - 60455-970 Fortaleza-Ce, Brazil.*

In this work, it is introduced a new function based on the non-trivial zeros of the Riemann-zeta function. Such function shows an interesting behavior: when the argument of the function grows, it changes from a pseudo-random behavior to a harmonic behavior with decreasing frequency.

## 1. Introduction

Let $s_1$, $s_2$, $s_3$, …, $s_k$, be a set of the first $k$ non-trivial zeros of the Riemann-zeta function, a function of a complex variable $s$ that analytically continues the sum of the infinite serie $\zeta(s) = \sum_n (1/n^s)$, where $s_j = 1/2 + ib_j$. Here, instead of using the zeros, we are going to use the variables $\theta_j$, where $e^{i\theta_j} = s_j^*/s_j$, for $j=1,…,k$ [1]. Hence

$$\theta_j = -\pi + \tan^{-1}\left[-b_j/\left(1/4 - b_j^2\right)\right]. \tag{1}$$

Now, one can define the function

$$y(x) = \frac{1}{x}\sum_{j=1}^{k} e^{ix(\pi+\theta_j)}. \tag{2}$$

An interesting behavior appears when one considers the function $\cos[y(x)]$. A plot of the real and imaginary parts of $\cos[y(x)]$ can be seen in Figs. 1 and 2 ($k = 169{,}165$).

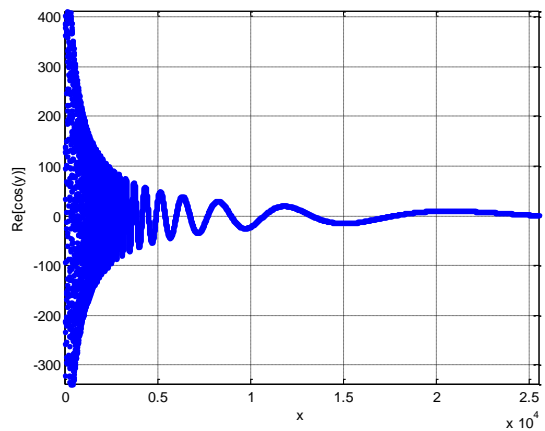

Fig. 1 – Real part of cos(*y*) versus *x*.

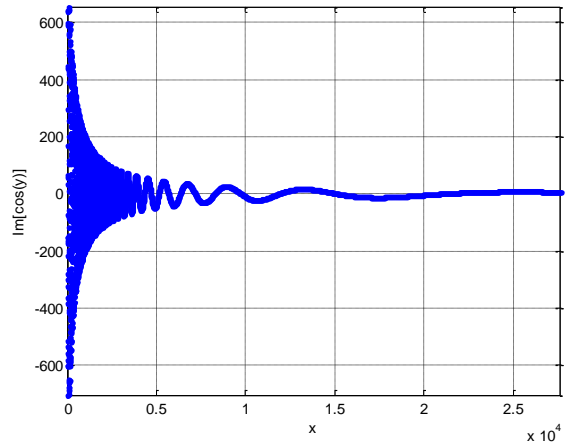

Fig. 2 – Imaginary part of cos(*y*) versus *x*.

The phase space is shown in Fig. 3.

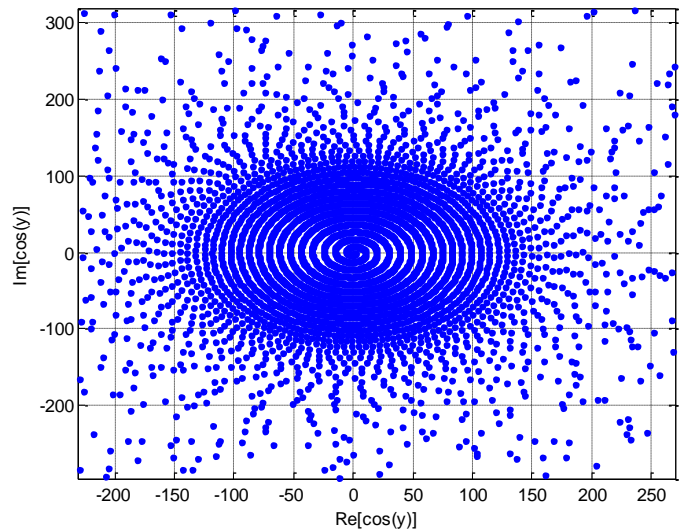

Fig. 3 – Re[cos(y)] versus Im[cos(*y*)].

As it can be noted in Figs.1, 2 and 3, for small values of *x* Re[cos(*y*)] and Im[cos(*y*)] show a pseudo-random behavior (the scattered dots at the phase space). When *x* grows, a harmonic behavior takes place (the spiral at the center of the phase space). Additionally, there is an amplitude damping and the frequency (in the harmonic region) decreases when *x* increases. A closer look at the pseudo-random and harmonic regions can be seen in Figs. 4-9 below.

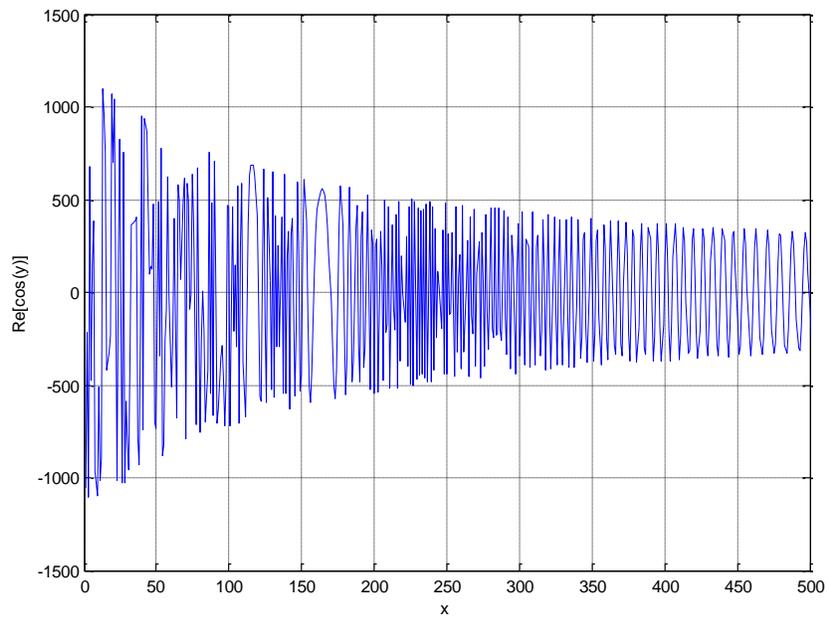

Fig. 4 – Re[cos(y)] versus $x \in [1,500]$: pseudo-random behavior.

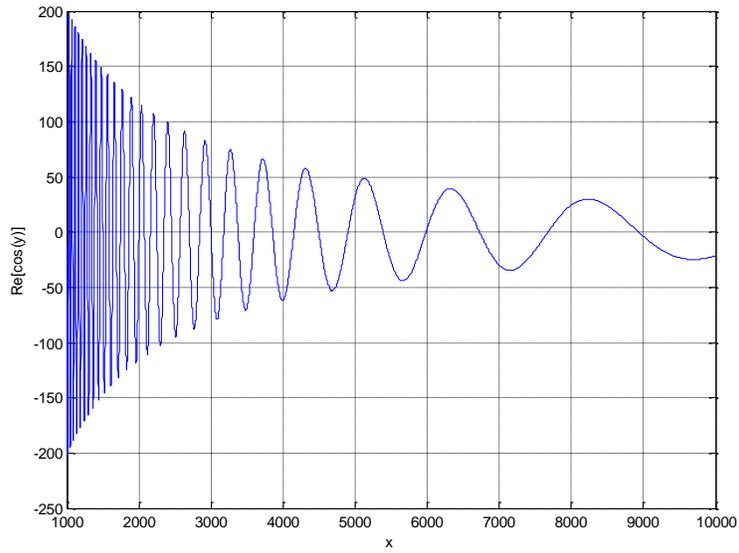

Fig. 5 – Re[cos(y)] versus $x \in [1000,10000]$: harmonic behavior with decreasing frequency.

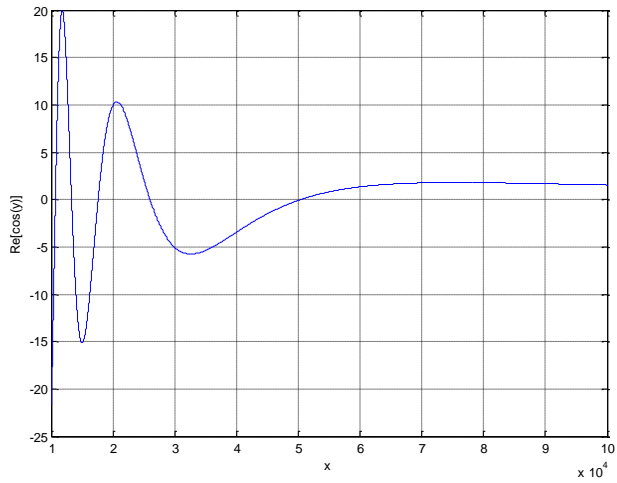

Fig. 6 – Re[cos(y)] versus $x \in [10000, 100000]$: harmonic behavior with decreasing frequency.

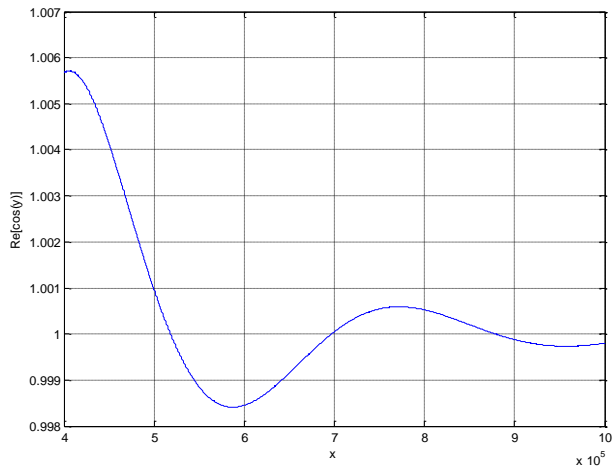

Fig. 7 – Re[cos(y)] versus $x \in [400000, 1000000]$: harmonic behavior with decreasing frequency.

At last, a 3D plot and a plot of the phase space where the harmonic part is clearly seen (the spiral) are shown in Figs. 8 and 9. At the phase space, the spiral tends to the point (0,0) because of the amplitude damping.

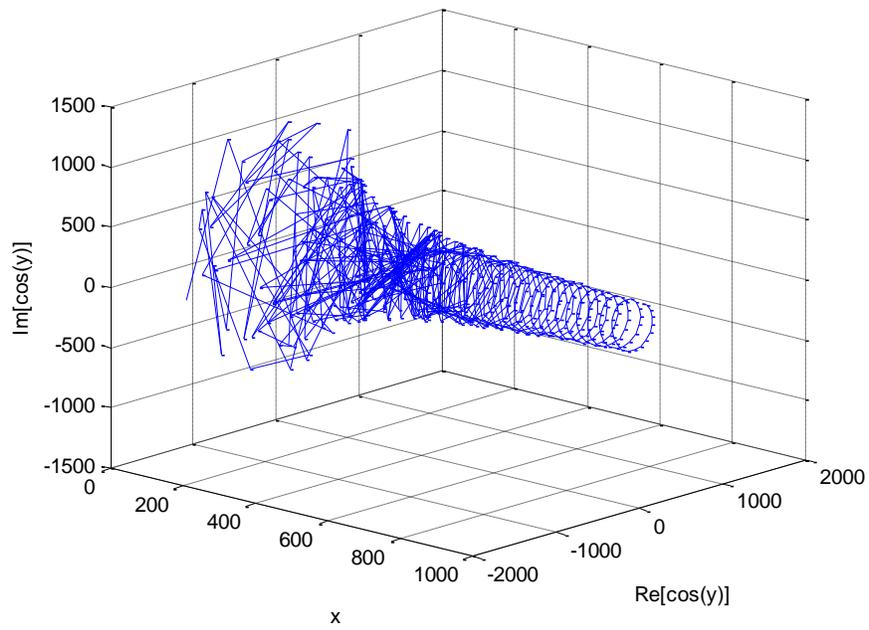

Fig. 8 – Re[cos(y)] versus Im[cos(y)] versus $x \in [1,1000]$: From pseudo-random to harmonic behavior.

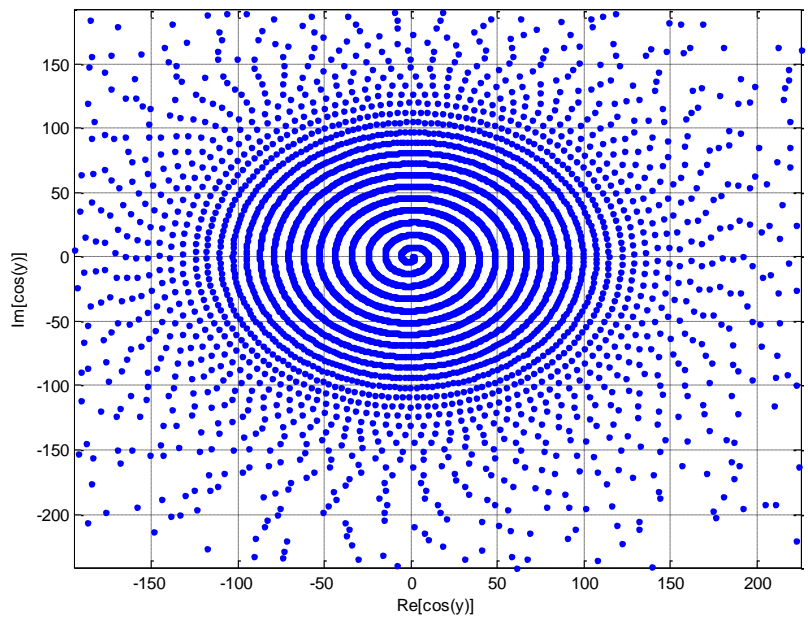

Fig. 9 – Re[cos(y)] versus Im[cos(y)]: The spiral represents the harmonic behavior and the scattered dots represents the pseudo-random region.

## 2. Conclusions

The function Re[cos($y(x)$)], with the argument $x$ playing the role of time, seems to be appropriate to model a non-ideal damped oscillator with frequency modulation, since it presents a transient (the pseudo-random region) followed by an oscillatory behavior with decreasing amplitude and frequency. However, it is not a trivial task to find out what physical system can be represented by Re[cos($y(x)$)] since: 1) It is not obvious to define, in a deterministic physical systems, what can generate the transient that looks like a pseudo-random signal. 2) the amplitude decay is not exponential; 3) the frequency variation depends on the zeros of the Riemann-zeta function. If such physical system exists, then one would have a new application of the Riemann-zeta function in physics [2] and engineering. In particular, if at least part of the frequency (in Fig. 5, for example) falls inside the region of human audibility (vision), one could hear (see) 'the sound of the zeros' ('the colors of the zeros').

Obviously, there are a lot of open questions. Two of them are: 1) What happens if more zeros are used? 2) How the frequency variation depends on the zeros?